\documentclass[11pt, reqno]{amsart}
\usepackage{amsmath, amsthm, amscd, amsfonts, amssymb, graphicx, color, upgreek}
\usepackage[bookmarksnumbered, plainpages]{hyperref}
\usepackage{enumerate , enumitem}

\usepackage[left=2.5 cm, right= 2.5 cm, top = 3 cm, bottom = 2.5 cm]{geometry}
\usepackage{parskip}

\def\R{\mathbb{R}}
\def\eps{\varepsilon}

\def\N{\mathbb{N}}

\def\S{\mathbb{S}}
\def\lt{\left}
\def\rt{\right}

\newtheorem{theorem}{Theorem}[section]
\newtheorem{lemma}[theorem]{Lemma}

\theoremstyle{definition}

\newtheorem{definition}[theorem]{Definition}

\theoremstyle{remark}
\newtheorem{remark}[theorem]{Remark}
\numberwithin{equation}{section}

\begin{document}
\setcounter{page}{1}

\centerline{}

\centerline{}


\title[A. Prade]{Boundary regularity for nonlocal elliptic equations \\ over Reifenberg flat domains}

\author[A. Prade]{Adriano Prade}

\address{CMAP, \'Ecole Polytechnique, Institut Polytechnique de Paris, 91120 Palaiseau, France.}
\email{adriano.prade@polytechnique.edu}

\begin{abstract} We prove sharp boundary regularity of solutions to nonlocal elliptic equations arising from operators comparable to the fractional Laplacian over Reifenberg flat sets and with null exterior condition. More precisely, if the operator has order $2s$ then the solution is $C^{s-\eps}$ regular for all $\eps>0$ provided the flatness parameter is small enough. The proof relies on an induction argument and its main ingredients are the construction of a suitable barrier and the comparison principle.
\end{abstract} 

\maketitle

\begin{section}{Introduction}

For $s \in (0,1)$ and $0 < \lambda \leq \Lambda < + \infty$, we consider linear nonlocal operators $L$ acting on functions $u: \R^n \rightarrow \R$ of the form
\begin{equation}\label{operators} \begin{split}
    Lu(x) & = \frac{1}{2} \int_{\mathbb{R}^n} \lt( 2u(x)-u(x+y)-u(x-y) \rt)K(y)\,dy \\ & = \text{P.V.} \int_{\mathbb{R}^n}(u(x)-u(x+y))K(y)\,dy,
\end{split}    
\end{equation}
with kernels $K: \mathbb{R}^n \rightarrow [0, +\infty]$ satisfying symmetry and strong ellipticity conditions, namely
    \begin{equation}\label{sym}
        K(y)= K(-y) \quad \text{and} \quad 0< \frac{\lambda}{|y|^{n+2s}} \leq K(y) \leq \frac{\Lambda}{|y|^{n+2s}}  \qquad \text{ for all } y \in \R^n.
    \end{equation}
Integro-differential operators of the type \eqref{operators}-\eqref{sym} are often known as operators comparable to the fractional Laplacian and their study has gained much attention in recent years. Indeed, nonlocal partial differential equations appear in multiple areas of mathematics and applied sciences such as analysis, probability, physics and biology, see \cite{fernandez2024integro, bucur2016nonlocal, garofalo2017fractional} and the examples presented therein. Among the many questions which have been investigated so far, we focus here on boundary regularity of solutions to nonlocal equations of the form
\begin{equation}\label{example} \begin{cases}
        Lu = f & \text{on } \Omega \\ 
        u=0 & \text{on } \R^n \setminus \Omega,
    \end{cases} \end{equation} 
where $L$ satisfies \eqref{operators}-\eqref{sym}, $\Omega \subset \R^n$ is a bounded open set with various additional assumptions and $f \in L^{\infty}(\Omega)$. It was shown in \cite{caffarelli2011regularity} that when $\Omega$ satisfies the exterior ball condition (so for example when $\Omega$ is a $C^{1,1}$ domain) then $u\in C^{\eps}(\overline{\Omega})$ for some $\eps>0$ small enough. The same was achieved in \cite[Prop. 2.6.9]{fernandez2024integro} for Lipschitz domains. Under the additional 
 hypothesis of homogeneity of the kernel
\[ K(y)= \frac{1}{|y|^{n+2s}} K\lt( \frac{y}{|y|}\rt),\]
sharp $C^s$ regularity up to the boundary of solutions to \eqref{example} was proven in \cite{fernandez2024integro, grube24} when $\Omega$ is respectively a $C^{1,\alpha}$ domain or satisfies a $C^{1,\text{Dini}}$ condition. This result for $C^{1,\alpha}$ domains was extended to operators comparable to the fractional Laplacian in \cite{ros2024optimal}, where it was also shown that $u\in C^{s-\eps}(\overline{\Omega})$ for all $\eps>0$ when $\Omega$ is $\delta$-Lipschitz for sufficiently small $\delta$. Many other regularity results under different assumptions are available in the literature, see \cite{grube24, ros2024optimal} for more detailed overviews.

The $C^{s-\eps}$ boundary regularity for $\delta$-Lipschitz sets from \cite{ros2024optimal} proves to be optimal and the purpose of this note is to generalize it to an even less smooth framework. In particular, we deal with the class of Reifenberg flat domains according to the following definition.
\begin{definition}\label{reif}
Let $ 0< \eta < 1/2 $ and $r_0 > 0$. A set $\Omega \subset \mathbb{R}^n$ is $(\eta, r_0)$-Reifenberg flat if for all $x \in \partial \Omega$ and for all $0<r<r_0$ there exists a hyperplane $H_{x,r}$ containing $x$ and such that:
    \begin{itemize}
        \item $d_H ( \partial \Omega \cap B_r(x), H_{x,r} \cap B_r(x) ) \leq \eta r$, where $d_H$ denotes the Hausdorff distance;
        \item one of the connected components of $\{d(\cdot, H_{x,r}) \geq 2\eta r \} \cap B_{r}(x)$ is included in $\Omega$ and the other in $\Omega^c$.
    \end{itemize}
    We say that $\Omega \subset \R^n$ is a Reifenberg flat domain whenever it is a bounded, open and connected Reifenberg flat set.
\end{definition}
The boundary of a Reifenberg flat set can be seen as a perturbation of a hyperplane at small scales, so their regularity may be rather mild. For example, the class includes both $\delta$-Lipschitz sets with sufficiently small $\delta$ and more general sets whose boundary has a fractal behaviour. Reifenberg flat sets have had a valuable role in many areas of analysis such as minimal surfaces, harmonic analysis and the study of elliptic and parabolic partial differential equations, see \cite{lemenant2012extension, lemenant2013boundary} for many references pertaining these fields. In this work we aim at proving the following.
\begin{theorem}\label{goal}
    Let $L$ be of the type \eqref{operators}-\eqref{sym} and $\eps \in (0,s)$. Let $\Omega \subset \R^n$ be a $(\eta, r_0)$-Reifenberg flat domain for some $r_0>0$, $f \in L^{\infty}(\Omega)$ and $u$ be the weak solution to
    \begin{equation}\label{solut} \begin{cases}
        Lu = f & \text{on } \Omega \\ 
        u=0 & \text{on } \R^n \setminus \Omega.
    \end{cases} \end{equation}
    Then, there exist $\eta_0 = \eta_0(n,s,\lambda, \Lambda, \eps)>0$ and $C= C(n,s, \lambda, \Lambda, \eps)>0$ such that for all $\eta \leq \eta_0$ we have $u \in C^{s-\eps}(\overline{\Omega})$ with
   \begin{equation}\label{weakbound} \|u\|_{C^{s - \eps}(\overline{\Omega})} \leq C \| f\|_{L^{\infty}(\Omega)}. \end{equation}
\end{theorem}
Many techniques have been employed so far to prove boundary regularity of elliptic equations over $(\eta, r_0)$-Reifenberg flat domains with $\eta$ small enough, for both local and nonlocal operators. Regarding the local case, boundary regularity for the Poisson equation with $L^q$ datum and homogeneous boundary condition was achieved in \cite{lemenant2013boundary}, by relying on a variant of the Alt-Caffarelli-Friedman monotonicity formula and on Campanato criterion. The result was recently extended to elliptic operators of order $2m$ with constant coefficients in \cite{lemenant2025} using variational methods, see also \cite{dong2011higher} for another generalization involving weaker hypotheses on the coefficients. Otherwise, boundary regularity is often proved by introducing a suitable barrier function, which bounds the growth of a solution close to the boundary via comparison principle. This was effectively carried out in \cite{lian2018boundary}, where various versions of the ABP maximum principle are applied to prove H\"older continuity for many classes of linear and nonlinear elliptic equations. In the same work, an iterative argument exploiting the flatness property was used in order to conclude.

Barrier methods prove successful in the nonlocal case as well, see \cite{fernandez2024integro} for their application to nonlocal equations over $C^{1,\alpha}$, $C^1$ and Lipschitz domains. Concerning Reifenberg flat sets, a first result in the spirit of the previous ones was achieved in \cite{goldman2024rigidity} by adapting the techniques of \cite{lemenant2013boundary}. The authors proved boundary regularity of the $1$-potential function $u_{\Omega}$ of a Reifenberg flat domain $\Omega$, which is a solution to the equation $(-\Delta)^{1/2}u=0$ in $\Omega^c$. Given $s\in (0,1)$ and $f\in C({\overline{\Omega}})$, the same conclusion was obtained in \cite{prade2024abp} for solutions to 
\[ \begin{cases} (-\Delta)^s u = f \qquad & \text{on } \Omega \\ u=0 & \text{on } \R^n \setminus \Omega, \end{cases} \]
by following the ideas of \cite{lian2018boundary}, thanks to a $2s$-potential function employed as barrier and a nonlocal version of the ABP maximum principle. Despite being a first step towards boundary regularity of nonlocal operators over Reifenberg flat domains, the result does not prove fully satisfactory. Indeed, it is limited to the specific case of the fractional Laplacian $(-\Delta)^s$ with uniformly continuous right-hand side, thus it is not much in accordance with the well-established theory for $C^{1,\alpha}$, $C^1$ and Lipschitz domains presented in \cite{fernandez2024integro}. Hence, Theorem \ref{goal} completes the picture initiated in \cite{prade2024abp} by extending boundary regularity over Reifenberg flat domains to operators of the form \eqref{operators}-\eqref{sym} with $L^{\infty}$ right-hand side. Before providing an outline of the proof, we also mention \cite{korvenpaa2016obstacle} for another H\"older continuity result for supersolutions to some linear equations, valid in a non-smooth setting including Reifenberg flat domains as well.

The proof of Theorem \ref{goal} follows the same approach adopted in \cite{ros2024optimal} for $\delta$-Lipschitz domains $\Omega$ with sufficiently small $\delta$. The authors show first $C^{\sigma}$ regularity for some $\sigma>0$ by introducing a barrier function and applying the maximum principle as mentioned before, see \cite[Prop. 2.6.4 and Prop. 2.6.9]{fernandez2024integro} for more details. Thanks to the geometry of the domain the construction is global, in the sense that the requested comparison between the solution and the barrier is obtained at once for any $x \in \Omega$ regardless of its distance from the boundary $\partial \Omega$. However, even though the barriers employed there are available also for Reifenberg flat domains (see \cite[Lemma B.3.3]{fernandez2024integro}), it is not clear how to adapt the same strategy in order to prove boundary regularity in this case as well. To overcome this issue, in Theorem \ref{th} we rely on an induction argument following the ideas of \cite{lian2018boundary, prade2024abp}. More precisely, we introduce first a barrier with respect to a flat boundary by exploiting the supersolutions built in Lemma \ref{finalbarrier}. Then, proceeding by induction allows us to get information on the behaviour of the solution at every scale by means of the weak comparison principle. In this way $C^{\sigma}$ boundary regularity of the solution to \eqref{solut} is established for some $\sigma >0$ and we also obtain $u \in C^{\sigma}(\overline{\Omega})$ by combining it with interior estimates. Finally, in Theorem \ref{compactness} we argue by compactness like in \cite{ros2024optimal} and we upgrade to $C^{s-\eps}$ regularity for any $\eps>0$.

The note is organized as follows. In Section \ref{section2} we set the notation and we list some preliminary results. In Section \ref{section3} we build supersolutions in quite general domains of $\R^n$ for our class of operators. In Section \ref{section4} we prove Theorems \ref{th} and \ref{compactness} thus getting boundary H\"older regularity. 
    
\end{section}

\begin{section}{Preliminaries}\label{section2}

In this section we gather all the definitions and preliminary results regarding nonlocal elliptic equations which will be useful throughout the note, the main reference is the monograph 
\cite{fernandez2024integro}.

We work in $\R^n$ with $n \geq 2$, we write $x = (x_1, ... , x_n)$ for $x \in \R^n$ and $\{e_i\}_{i=1}^n$ to denote the standard orthonormal basis. Expression $\chi_A$ stands for the indicator function of $A\subset \R^n$, if $A$ is measurable then $|A|$ denotes its Lebesgue measure. For $x\in \R^n$ and $r>0$, $B_r(x):=\{y\in \R^n: |x-y|<r\}$. We drop the dependence on $x$ if $x=0$ and we set $\omega_n:=|B_1|$. We often use constants $C>0$ or $c>0$ for estimates, their dependence is given in the corresponding statement or proof and their value may change form line to line. We also denote $A \ll B$ to indicate that there exists a (typically small) universal constant $\eps>0$ depending only on $n$ such that $A \leq \eps B$. Finally, we use the convention that if $\beta=k+\alpha$ for $k \in \N$ and $\alpha\in(0,1)$ then $C^{\beta}(\Omega)$ coincides with the H\"older space $C^{k,\alpha}(\Omega)$.

For $\Omega \subset \R^n$ open and $x \in \Omega$, let $d_{\Omega}(x) := \text{dist}(x, \Omega^c)$. By \cite[Lemma B.0.1]{fernandez2024integro}, there exist a function $\updelta_{\Omega} \in C^{\infty}(\Omega) \cap \text{Lip}(\R^n)$ and constants $C, C_k>0$ depending only on $n$ and $k \in \N$ such that
\begin{equation} \label{regdist}
    d_{\Omega} \leq \updelta_{\Omega} \leq C d_{\Omega} \qquad \text{on }\, \R^n \end{equation}
and
\begin{equation} \label{regdist1} \lt \| \updelta_{\Omega} \rt \|_{C^k(\Omega)} \leq C_k d_{\Omega}^{1-k} \qquad \text{for all } k \in \N. \end{equation}
When it is clear from the context, we drop the dependence on $\Omega$ by simply writing $d$ or $\updelta$. 

As anticipated, for $s \in (0,1)$ and $0 < \lambda \leq \Lambda < + \infty$ we deal we the class of nonlocal elliptic operators
\[ \mathcal{L}_s^n(\lambda, \Lambda) := \lt \{L \text{ of the form \eqref{operators} and with kernel \textit{K} satisfying \eqref{sym}} \rt \}. \]
The class $\mathcal{L}_s^n(\lambda, \Lambda)$ is scale invariant. Namely, if $L \in \mathcal{L}_s^n(\lambda, \Lambda)$ and $u_r(x):=u(rx)$ for $r>0$ then there exist $L_r \in \mathcal{L}_s^n(\lambda, \Lambda)$ such that
\begin{equation}\label{rescaledop}
    (L_ru_r)(x)=r^{2s}(L u)(rx).
\end{equation}
If $L$ has kernel $K$, a change of variables shows that $L_r$ has kernel $K_r(y) := r^{n+2s}K(ry)$. This also shows that $L$ is scale invariant if and only if its kernel is homogeneous of degree $-(n+2s)$. Now, we introduce the function space
\[ L^1_{2s}(\R^n):= \lt \{ u: \R^n \longrightarrow \R : \|u\|_{L^1_{2s}(\R^n)}:= \int_{\R^n} \frac{|u(y)|}{1+|y|^{n+2s}}dy < + \infty\rt \}\]
and we provide a sufficient condition in order to evaluate $Lu$ in the $L^{\infty}$ sense, see Lemma 2.2.4 and Lemma 1.2.3 from \cite{fernandez2024integro}.
\begin{lemma}\label{wlog}
    Let $L \in \mathcal{L}_s^n(\lambda, \Lambda)$ and $u\in C^{2s+\eps}(B_1) \cap L^1_{2s}(\R^n)$ for some $\eps>0$, then $Lu\in L^{\infty}_{loc}(B_1)$ and there exists $C=C(n,s,\eps)>0$ such that
    \[ \| Lu \|_{L^{\infty}}(B_{1/2}) \leq C\Lambda \lt( \| u\|_{C^{2s+\eps}(B_1)} + \| u\|_{L^1_{2s}(\R^n)}\rt).\]
\end{lemma}
Next, we present two different notions of solution to nonlocal elliptic equations. We recall first the definition of the fractional Sobolev space: 
\[ H^s(\R^n):= \lt \{ u\in L^2(\R^n): [u]_{H^s(\R^n)}:= \lt( \int_{\R^n \times \R^n} \frac{(u(x)-u(y))^2}{|x-y|^{n+2s}}dx\,dy \rt)^{1/2}< + \infty\rt \}. \]
For any kernel $K$ satisfying \eqref{sym} and $\Omega \subset \R^n$ it is well-defined the symmetric and positive semi-definite bilinear form $\langle \cdot, \cdot \rangle_{K,\Omega} : C^{\infty}_c(\mathbb{R}^n)\times C^{\infty}_c(\mathbb{R}^n) \longrightarrow \mathbb{R}$ given by
\[ \langle u, v \rangle_{K,\Omega} := \frac{1}{2} \int_{\mathbb{R}^n \times \mathbb{R}^n\setminus (\Omega^c\times \Omega^c)} \left( u(x)-u(y)\right) \left( v(x)-v(y) \right) K(x-y)\,dx\,dy. \]
Thanks to this formula, a linear nonlocal operator of the form \eqref{operators} is well-defined via the relation $(Lu,\varphi):=\langle u, \varphi \rangle_{K,\Omega}$ for $\varphi \in H^s(\R^n)$ with $\varphi=0$ on $\Omega^c$, so it is natural to give the following definition of weak solution.
\begin{definition}[Weak solution]
    Let $ L \in \mathcal{L}_s^n(\lambda, \Lambda)$ and let $\Omega \subset \R^n$ be a bounded open set. Let $f\in L^{\infty}(\Omega)$, then a function $u$ such that $ \langle u, u\rangle_{K,\Omega} < + \infty$ is a weak subsolution to $Lu \leq f$ in $\Omega$ if
    \begin{equation}\label{weaksol}
        \langle u, \varphi \rangle_{K,\Omega}  \leq \int_{\Omega}f\varphi \quad \text{for all } \varphi \in 	H^s(\mathbb{R}^n) \, \text{ with } \varphi\geq 0 \text{ and }\varphi =0 \mbox{ on } \mathbb{R}^n \setminus \Omega.
    \end{equation}
    We say that $u$ is a weak supersolution to $Lu \geq f$ in $\Omega$ if \eqref{weaksol} holds with $\geq$ instead of $\leq$. Moreover, $u$ is a weak solution to $Lu = f$ in $\Omega$ if it is both a weak subsolution and a weak supersolution.
\end{definition}
See \cite[Theorem 2.2.24]{fernandez2024integro} for an existence result and \cite[Lemma 2.3.9]{fernandez2024integro} for an $L^{\infty}$ bound for weak solutions to the Dirichlet problem. It is possible to prove that $u L\varphi \in L^1(\R^n)$ for all $\varphi \in C^{\infty}_c(\Omega)$ whenever $u\in L^1_{2s}(\R^n)$, so we can weaken the notion of solution even more.
\begin{definition}[Distributional solution]
    Let $ L \in \mathcal{L}_s^n(\lambda, \Lambda)$ and let $\Omega \subset \R^n$ be a bounded set. Let $f\in L^1_{\text{loc}}(\R^n)$, then $u\in L^1_{2s}(\R^n)$ is a distributional subsolution to $Lu \leq f$ in $\Omega$ if
    \begin{equation}\label{distsol} \int_{\R^n} uL\varphi \leq \int_{\Omega} f\varphi \qquad \text{for all } \varphi \in C^{\infty}_c(\Omega) \text{ with } \varphi \geq 0.\end{equation}
    We say that $u$ is a distributional supersolution to $Lu \geq f$ in $\Omega$ if \eqref{distsol} holds with $\geq$ instead of $\leq$. Moreover, $u$ is a distributional solution to $Lu = f$ in $\Omega$ if it is both a distributional subsolution and supersolution.
\end{definition}
If $f\in L^{\infty}(\R^n)$ and $\langle u, u\rangle_{K,\Omega} < + \infty$, then $Lu \leq f$ in $\Omega$ in the weak sense if and only if it is true also in the distributional sense. We now present a very powerful tool for studying nonlocal partial differential equation, the maximum principle for weak solutions.
\begin{theorem}[Weak maximum principle]\label{maxprinciple}
    Let $L \in \mathcal{L}_s^n(\lambda, \Lambda)$ and let $\Omega \subset \R^n$ be a bounded open set. Let $u$ be a weak solution to $Lu\geq 0$ in $\Omega$ with $u\geq 0$ in $\R^n \setminus \Omega$, then $u \geq 0$ in $\R^n$.
\end{theorem}
Its immediate consequence, the weak comparison principle, will be used in Section \ref{section4} to prove boundary regularity. We also point out that the same result holds for distributional solutions which are continuous up to the boundary, namely for $u\in C(\overline{\Omega}) \cap L^1_{2s}(\R^n)$, see \cite[Lemma 2.3.5]{fernandez2024integro}. Finally, we recall an interior regularity result for weak and distributional solutions. 
\begin{theorem}[Interior regularity]\label{interiorestimates}
    Let $L \in \mathcal{L}_s^n(\lambda, \Lambda)$ and $s \neq \frac{1}{2}$ . Let $f \in L^{\infty}(B_1)$ and $u \in L^1_{2s}(\R^n)$ be a distributional solution to $Lu=f$ in $B_1$, then $u \in C^{2s}_{\text{loc}}(B_1)$ and there exists $C=C(n, s, \lambda, \Lambda)>0$ such that
    \[ \|u\|_{C^{2s}(B_{1/2})} \leq C \lt( \| u \|_{L^1_{2s}(\R^n)} + \,\| f\|_{L^{\infty}(B_1)} \rt). \]
    If instead $s=\frac{1}{2}$, then $u \in C^{1-\eps}(B_1)$ for all $\eps>0$ and the same estimate holds with $\|u\|_{C^{1-\eps}(B_{1/2})}$.
\end{theorem}
\end{section}

\begin{section}{Barriers}\label{section3}
The goal of this section is to build barriers for open sets $\Omega \subset \R^n$ satisfying the following property for some $R>0$ and $0<\kappa<R$. We will apply these results in Section \ref{section4} in the case when $\Omega$ is the complement of a ball.
\begin{enumerate}[label=$\lt( \mathcal{P}_{R,\kappa} \rt)$, ref=$\lt( \mathcal{P}_{R,\kappa} \rt)$]
    \item\label{prop} For all $0<r<R$ and $z \in \partial \Omega$ there exists $x_{r,z} \in \Omega^c$ for which $B_{\kappa r}(x_{r,z}) \subset \Omega^c \cap B_r(z)$.
\end{enumerate}
We point out that in assumption (B.3.1) from \cite[Appendix B.3]{fernandez2024integro} the same condition is asked but for all $r>0$ instead of only $0<r<R$. Here we adapt to this weaker framework some of the technical results presented there, starting with the following GMT lemma.
\begin{lemma}\label{GMT}
Let $\Omega \subset \R^n$ be an open set satisfying \ref{prop} for some $R>0$ and $0<\kappa<R$. Then, there exist $\theta = \theta(n,\kappa) \in (0,1)$ and $C=C(n,\kappa)>0$ such that
\[ \big| \{ 0 < d_{\Omega} < r\rho \} \cap B_{\rho}(z)\big| \leq C r^{\theta} |B_{\rho}| \]
for all $z \in \partial \Omega$, $r \in (0,1/2)$ and $\rho \in (0,2R)$.
\end{lemma}

\begin{proof} Given $z \in \partial \Omega$ and $\rho \in (0,2R)$, by considering the set $\Omega_{\rho} := \frac{\Omega - z}{\rho}$ it is enough to prove that 
\[ \big| \{ 0 < d_{\Omega_{\rho}} < r \} \cap B_1 \big| \leq C r^{\theta} \qquad \text{for all } r \in (0,1/2) \]
for some $\theta \in (0,1)$ and $C>0$ depending on $n$ and $\kappa$. Noticing that $\Omega_{\rho}$ satisfies $(P_{1/2,\kappa})$, we argue as in the proof of \cite[Lemma B.3.5]{fernandez2024integro} and we conclude. \end{proof}

The next result is the equivalent of \cite[Lemma B.3.2]{fernandez2024integro} for open sets $\Omega$ satisfying \ref{prop}. The same uniform estimate holds but only for distances up to order 1 with respect to either $R$ or the diameter of $\Omega^c$.

\begin{lemma}\label{lemmabarrier} Let $ L \in \mathcal{L}_s^n(\lambda, \Lambda)$ and let $\Omega \subset \R^n$ be an open set satisfying \ref{prop} for some $R>0$ and $0<\kappa<R$. Then, there exist $d_0 = d_0(s, R, \Omega) >0$ and $c_0 = c_0 (n, s, \lambda, \kappa)>0$ such that
\[ L\chi_{\Omega}(x) \geq c_0 \,d^{-2s}(x) \qquad \textit{for all } x \in \Omega \textit{ for which } d(x) < d_0. \]
\end{lemma}

\begin{proof}
    Let $D:=\text{diam}(\Omega^c)$. If $D < + \infty$, up to decreasing $\kappa$ we may assume $R \gg D$ without loss of generality. We set
    \[ d_0 := \min \lt \{ 10D, \frac{R}{2^{1/2s+1}} \rt \}.\]
    Let $x \in \Omega$ be such that $d_{\Omega}(x) < d_0$, we assume $x=0$ and denote $r_0:=d_{\Omega}(0)$. Take $z \in \partial \Omega$ for which $r_0 = |z|$, then $B_{r-r_0}(z) \subset B_r$ for all $r > r_0$ and \ref{prop} implies that
    \begin{equation}\label{ballestimate}
    \big| B_r \cap \Omega^c \big| \geq \big| B_{r-r_0}(z) \cap \Omega^c \big| \geq \kappa^n |B_{r-r_0}| = \kappa^n \lt( \frac{r-r_0}{r}\rt)^n |B_r| \qquad \text{for } r_0<r<R.
    \end{equation} 
    We compute
    \[ \begin{split}
        L\chi_{\Omega}(0) & = \text{P.V.} \int_{\R^n} \big( \chi_{\Omega}(0) - \chi_{\Omega}(y) \big) K(y)dy = \int_{\R^n} \chi_{\Omega^c}(y)K(y)dy \geq \lambda \int_{\Omega^c} \frac{dy}{|y|^{n+2s}}.
    \end{split}  \]
    By layer cake representation and a change of variables, we find
    \[ \begin{split}
        \int_{\Omega^c} \frac{dy}{|y|^{n+2s}} & = \int_0^{+ \infty} \lt| \Omega^c \cap \lt \{ \frac{1}{|y|^{n+2s}} > t \rt\}\rt| dt = \int_0^{+ \infty} \lt| \Omega^c \cap B_{t^{-1/(n+2s)}}\rt| dt \\ & = (n+2s) \int_0^{+ \infty}  \lt|\Omega^c \cap B_{\rho} \rt| \frac{d\rho}{\rho^{n+2s+1}} \geq \int_{r_0}^{+\infty} \lt| \Omega^c \cap B_{\rho} \rt| \frac{d\rho}{\rho^{n+2s+1}}.
    \end{split} \]
    If $r_0 \leq D/10$, from \eqref{ballestimate} we get
    \begin{equation}\label{Dinf1} |\Omega^c \cap B_{\rho}| \geq \frac{\kappa^n}{2^n} |B_{\rho}| = c_1 \rho^n \qquad \text{for all } 2r_0 < \rho < R.  \end{equation}
    Altogether,
    \begin{equation}\label{Dinf2}  L\chi_{\Omega}(0) \geq \lambda \int_{2r_0}^R \lt| \Omega^c \cap B_{\rho} \rt| \frac{d\rho}{\rho^{n+2s+1}} \geq c_1 \int_{2r_0}^R \frac{d\rho}{\rho^{2s+1}} \geq c_1 \lt( r_0^{-2s} - 2^{2s}R^{-2s} \rt) \geq c_1 \, r_0^{-2s}. \end{equation}
    Otherwise, if $D/10 < r_0 < 10D$ then by \eqref{ballestimate} we obtain
    \[ |\Omega^c \cap B_{\rho}| \geq \frac{\kappa^n}{1000^n} |B_{\rho}| = c_2 \rho^n \qquad \text{for all } \quad r_0 + D/10 < \rho < 100D. \]
    Similarly to before, we have
    \[ L\chi_{\Omega}(0) \geq \lambda \int_{r_0  + D/10}^{100D} \lt| \Omega^c \cap B_{\rho} \rt| \frac{d\rho}{\rho^{n+2s+1}} \geq c_2 \int_{2r_0}^{10r_0} \frac{d\rho}{\rho^{2s+1}} \geq c_2 \lt( r_0^{-2s} - 5^{-2s}R^{-2s} \rt)  \geq c_2 \, r_0^{-2s}. \]
    We finish the proof in the case $D<+\infty$ by setting $c_0 = c_0(n, s, \lambda, \kappa) := \max\{c_1, c_2\}$. If instead $D=+\infty$ then we always have $r_0 \leq D/10$, so we reach the same conclusion by arguing like in \eqref{Dinf1} and \eqref{Dinf2}.
\end{proof}
\begin{remark}\label{remarkdistance}
    A quick inspection of the proof shows that if $c_0$ is allowed to vary depending on $d(x)$ then the conclusion $L\chi_{\Omega}(x) \geq c_0 \,d^{-2s}(x)$ remains valid for any $x \in \Omega$. More precisely, up to changing $d_0$ and adequately modifying the ensuing computations we see that $c_0 \rightarrow 0$ as $d_0 \rightarrow +\infty$. The precise choice of $d_0$ made above will be clear in section \ref{section4}, applying Lemma \ref{finalbarrier} when $\Omega$ is the complement of a ball.
\end{remark}
Thanks to the previous lemma we can construct supersolutions up to distances of order 1 in quite general domains. Once Lemma \ref{GMT} and Lemma \ref{lemmabarrier} are established the proof of the next result is very similar to \cite[Lemma B.3.3]{fernandez2024integro} and we reproduce it for convenience of the reader. Lemma \ref{finalbarrier} will be useful later in Section \ref{section4} for proving Theorem \ref{th}.
\begin{lemma}\label{finalbarrier} Let $ L \in \mathcal{L}_s^n(\lambda, \Lambda)$ and let $\Omega \subset \R^n$ be an open set with compact boundary satisfying \ref{prop} for some $R>0$ and $0<\kappa<R$. Then, there exist constants $\eps_0 = \eps_0(n, s, \lambda, \Lambda, \kappa)>0$ and $c=c(n, s, \lambda, \Lambda, \kappa)>0$ such that for all $\eps < \eps_0$
\begin{equation}\label{formula}
    L\updelta^{\eps}(x) \geq c\, d^{\eps - 2s}(x) \qquad \textit{for all } x \in \Omega \textit{ for which } d(x) < d_0,
\end{equation}
where $d_0$ is the same distance of Lemma \ref{lemmabarrier}.
\end{lemma}
\begin{proof}
    Since $\text{diam}(\partial \Omega)<+\infty$, as before we may also assume $R \gg \text{diam}(\partial \Omega)$ without loss of generality. Let $x \in \Omega$ such that $d(x) < d_0$, up to translation we write $x=0$ and set $2r_0:=d(0)$. Since $\updelta^{\eps}(0)\geq cr_0^{\eps}$, by Lemma \ref{lemmabarrier} we have
    \[ (L\updelta^{\eps})(0) = L(\updelta^{\eps}(0) \chi_{\Omega})(0) + L(\updelta^{\eps} - \updelta^{\eps}(0) \chi_{\Omega})(0) \geq c_0 \, r_0^{\eps - 2s} + L(\updelta^{\eps} - \updelta^{\eps}(0) \chi_{\Omega})(0). \]
    Hence, to obtain \eqref{formula} it is enough to show
    \begin{equation}\label{formula1} \big| L(\updelta^{\eps} - \updelta^{\eps}(0) \chi_{\Omega})(0)\big| \leq \frac{c_0}{2} r_0^{\eps - 2s}. \end{equation}
    By property \eqref{regdist1} of $\updelta$, a Taylor expansion shows that for $\eps>0$ small
    \[ \lt \| \updelta^{\eps} - \updelta^{\eps}(0) \chi_{\Omega} \rt \|_{C^{2}(B_{r_0})} \leq C \eps r_0^{\eps-2}.\]
    Therefore
    \[ \begin{split}
        \big| L(\updelta^{\eps} - \updelta^{\eps}(0) \chi_{\Omega})(0)\big| & \leq \Lambda \int_{B_{r_0}} \frac{C \eps r_0^{\eps-2}}{|y|^{n+2s-2}}dy + \Lambda \int_{\R^n \setminus B_{r_0}} \frac{\lt| (\updelta^{\eps} - \updelta^{\eps}(0) \chi_{\Omega})(y) \rt|}{|y|^{n+2s}}dy \\ & \leq C \eps r_0^{\eps - 2s} + \Lambda \int_{\Omega \setminus B_{r_0}} \frac{\lt| (\updelta^{\eps} - \updelta^{\eps}(0) \chi_{\Omega})(y) \rt|}{|y|^{n+2s}}dy,
    \end{split} \]
    where in the last line we used that $(\updelta^{\eps} - \updelta^{\eps}(0) \chi_{\Omega})(y)=0$ for all $y \in \Omega^c$. To bound the last integral we introduce the parameters $\eta, \gamma >0$ small and $M>1$ large to be chosen later. Then, we partition $\Omega \setminus B_{r_0}$ into
    \begin{itemize}
        \item[] $A_0 := \lt\{ y \in \R^n: |(\updelta^{\eps}(y) - \updelta^{\eps}(0))| \leq \eta\, \updelta^{\eps}(0) \rt\} \setminus B_{r_0}$ ,
        \item[] $D_0 := \lt\{ y \in \R^n: 0 < \updelta^{\eps}(y) < (1-\eta)\updelta^{\eps}(0) \rt\} \subset \lt \{ y\in \R^n: 0 < \updelta(y) < \gamma r_0  \rt\}$ ,
        \item[] $E_0 := \lt\{ y \in \R^n: \updelta^{\eps}(y) > (1+\eta)\updelta^{\eps}(0) \rt\} \subset (\Omega \setminus B_{Mr_0}$) .
    \end{itemize}
    Notice that the inclusions hold for $\eps >0$ small and depending only on $\eta$, $\gamma$ and $M$. Since $\updelta^{\eps}(0) \leq Cr_0^{\eps}$, we find
    \[ \int_{A_0} \frac{\lt| (\updelta^{\eps} - \updelta^{\eps}(0) \chi_{\Omega})(y) \rt|}{|y|^{n+2s}}dy \leq C\eta r_0^{\eps} \int_{\Omega \setminus B_{r_0}} \frac{dy}{|y|^{n+2s}} \leq C \eta r_0^{\eps-2s}. \]
    In addition,
    \[ \int_{E_0} \frac{(\updelta^{\eps} - \updelta^{\eps}(0) \chi_{\Omega})(y)}{|y|^{n+2s}}dy \leq \int_{\Omega \setminus B_{Mr_0}}\frac{\updelta^{\eps}(y)}{|y|^{n+2s}}dy \leq C\int_{\Omega \setminus B_{Mr_0}} \frac{1}{|y|^{n+2s-\eps}}dy \leq C (Mr_0)^{\eps-2s}. \]
    It remains to estimate the integral over $D_0$. As in the previous proof, by layer cake representation and a change of variables we obtain
    \[ \begin{split} \int_{D_0} \frac{(\updelta^{\eps}(0) \chi_{\Omega} - \updelta^{\eps})(y)}{|y|^{n+2s}}dy & \leq \int_{\lt \{0 < \updelta < \gamma r_0  \rt\}}\frac{\updelta^{\eps}(0)}{|y|^{n+2s}}dy \leq C r_0^{\eps} \int_{\lt \{0 < \updelta < \gamma r_0  \rt\}} \frac{dy}{|y|^{n+2s}} \\ & \leq C r_0^{\eps} \int_0^{+\infty} \lt| \lt \{0 < \updelta < \gamma r_0  \rt\} \cap B_{\rho} \rt| \frac{d\rho}{\rho^{n+2s+1}}. \end{split}\]
    Take $z \in \partial \Omega$ for which $2r_0 = |z|$, then $B_{\rho} \subset B_{3\rho}(z)$ for all $\rho > r_0$ and Lemma \ref{GMT} implies that there exist $\theta\in(0,1)$ and $C>0$ such that for every $r_0 < \rho < 2R$ we have
    \[ \lt|\lt \{0 < \updelta < \gamma r_0  \rt\} \cap B_{\rho} \rt| \leq \lt|\lt \{0 < d < \gamma r_0  \rt\} \cap B_{3\rho}(z) \rt| \leq  \lt|\lt \{0 < d < \gamma \rho  \rt\} \cap B_{3\rho}(z) \rt| \leq C\gamma^{\theta} \rho^n.  \]
    Recalling that $\partial \Omega \subset B_{2R}$ by assumption, we get
    \[ \int_{D_0} \frac{(\updelta^{\eps}(0) \chi_{\Omega} - \updelta^{\eps})(y)}{|y|^{n+2s}}dy \leq C \gamma^{\theta} r_0^{\eps} \int_{r_0}^{2R} \frac{d\rho}{\rho^{2s+1}} \leq C \gamma^{\theta} r_0^{\eps - 2s}. \]
    Thus, combining all the previous inequalities we obtain
    \[ \big| L(\updelta^{\eps} - \updelta^{\eps}(0) \chi_{\Omega})(0)\big| \leq C( \eps + \eta + M^{\eps-2s} + \gamma^{\theta})\,r_0^{\eps-2s}. \]
    We choose first $\eps$ small, then $\eta$, $\gamma$ small and $M$ large enough so that $C( \eps + \eta + M^{\eps-2s} + \gamma^{\theta})\leq c_0/2$, in this way \eqref{formula1} follows and the proof is concluded.
\end{proof}
\begin{remark}
    In view of Remark \ref{remarkdistance}, also the estimate $L\updelta^{\eps}(x) \geq c\, d^{\eps - 2s}(x)$ holds for any $x \in \Omega$ up to properly reducing the parameter $c$ when $d(x)$ increases.
\end{remark}
\end{section}

\begin{section}{Boundary regularity}\label{section4}

The purpose of this section is to prove $C^{s-\eps}$ regularity of solutions to nonlocal equations over $(\eta, r_0)$-Reifenberg flat sets with sufficiently small $\eta$. In Theorem \ref{th} we show first $C^{\sigma}$ regularity for a small $\sigma>0$. The proof goes by induction and its main ingredient is the application of the weak maximum principle Theorem \ref{maxprinciple}.

\begin{theorem}\label{th}
   Let $L \in \mathcal{L}_s^n(\lambda, \Lambda)$. Let $\Omega \subset \R^n$ be an open and connected $(\eta, r_0)$-Reifenberg flat set for some $r_0>0$, $f \in L^{\infty}(\Omega \cap B_2)$ and $u \in L^1_{2s}(\R^n)$ be a weak solution to
   \[ \begin{cases}
        Lu = f & \text{on } \Omega \cap B_2 \\ 
        u=0 & \text{on } B_2 \setminus \Omega.
    \end{cases} \]
   Then, there exist $\eta_0 = \eta_0(n,s,\lambda, \Lambda)>0$, $\sigma = \sigma(n,s,\lambda, \Lambda) >0$ and $C= C(n,s, \lambda, \Lambda)>0$ such that for all $\eta \leq \eta_0$ we have $u \in C^{\sigma}_{loc}(B_2)$ with
   \[ \|u\|_{C^{\sigma}(B_{1/2})} \leq C \lt( \| u \|_{L^1_{2s}(\R^n)} + \,\| f\|_{L^{\infty}(\Omega \cap B_2)} \rt). \]
\end{theorem}
Before starting the proof, we set some notation and we exploit the supersolutions we built in Section \ref{section3} to define suitable barriers for the problem. For $\eta>0$ small, consider the closed ball
\[ \mathcal{B}_0 := \overline{B \lt(- \lt(\frac{1}{8} +\eta \rt) e_n, \frac{1}{8} \rt)} \subset B_1. \]
Notice that $\mathcal{B}_0^c:=\R^n \setminus \mathcal{B}_0$ satisfies \ref{prop} for any $R>\kappa$ such that $R\kappa=1/8$, also $\text{diam}(\mathcal{B}_0)= 1/4$. Given some $\rho\in(0,1)$, we set $\mathcal{B}_k := \rho^k \mathcal{B}_0$ and $\mathcal{B}_k^c := \R^n \setminus \mathcal{B}_k$ for all $k \in \N$.  Let $K$ be the kernel of $ L \in \mathcal{L}_s^n(\lambda, \Lambda)$, for all $k\in \N$ we define the kernel $K_k$ as
\[ K_k(y) := \rho^{k(n+2s)}K(\rho^k y) \]
and we call $L_k$ the corresponding nonlocal operator. We have that $ L_k \in \mathcal{L}_s^n(\lambda, \Lambda)$ for $k\in \N$ by the scaling properties of the class. Hence, by Lemma \ref{finalbarrier} there exist $c, \eps_0 >0$ depending on $n$, $s$, $\lambda$ and $\Lambda$ such that for all $\eps \leq \eps_0$
\begin{equation}\label{barrier}
    L_k\updelta_{\mathcal{B}_0^c}^{\eps} \geq c \,d_{\mathcal{B}_0^c}^{\eps - 2s} \qquad \text{on } B_1  \qquad \text{for all } k\in \N.
\end{equation}
For $\eps < \eps_0$ to be chosen later, we take as barriers for our proof suitable rescalings of $\updelta_{\mathcal{B}_k^c}^{\eps}$: 
\[ v_k(x) := 4^{\eps}\rho^{(k-1)\sigma}\updelta_{\mathcal{B}_0^c}^{\eps}\lt( \frac{x}{\rho^k} \rt) \qquad \text{for } k\in\N. \]
Now, we gather some useful properties of the functions $v_k$. First of all, there holds
\begin{equation}\label{barrierk}
    Lv_k \geq c\, 4^{\eps}\rho^{(k-1)\sigma} \,d_{\mathcal{B}_k^c}^{\eps - 2s} \qquad \text{on } B_{\rho^k}  \qquad \text{for all } k\in \N.
\end{equation}
To show it, we recall that $d\lt( \tau x , E \rt) = \tau \, d \lt( x, \frac{1}{\tau} E \rt)$ for all $\tau > 0$ and $E \subset \R^n$, so for $\tau = \frac{1}{\rho^k}$ and $E = \mathcal{B}_0^c$ we get \begin{equation}\label{distscale} d_{\mathcal{B}_0^c}\lt( \frac{x}{\rho^k} \rt)  = \frac{1}{\rho^k} \, d_{\mathcal{B}_k^c}(x). \end{equation}
Thanks to the scaling property \eqref{rescaledop}, \eqref{barrier} and \eqref{distscale}, for $x\in B_{\rho^k}$ we obtain \[ \begin{split}
    L \lt( \updelta_{\mathcal{B}_0^c}^{\eps} \lt( \frac{x}{\rho^k} \rt) \rt) & = \lt( \frac{1}{\rho^k} \rt)^{2s} (L_k \updelta_{\mathcal{B}_0^c}^{\eps}) \lt( \frac{x}{\rho^k} \rt) \geq c\lt( \frac{1}{\rho^k} \rt)^{2s}  d_{\mathcal{B}_0^c}^{\eps-2s}\lt( \frac{x}{\rho^k} \rt) = c\lt( \frac{1}{\rho^k} \rt)^{\eps} d_{\mathcal{B}_k^c}^{\eps - 2s}(x) \\ & \geq c \,d_{\mathcal{B}_k^c}^{\eps - 2s}(x),
\end{split} \]
which yields \eqref{barrierk}. Next, by \eqref{regdist} there exists $C_H$ such that $\updelta_{\mathcal{B}_0^c}^{\eps} \leq C_H \eta^{\eps}$ when $d_{\mathcal{B}_0^c}(x) \leq 2\eta$. So, again by the scaling relation \eqref{distscale} we obtain
\begin{equation}\label{holderbarrier}
        v_k(x) \leq 4^\eps C_H \rho^{(k-1)\sigma} \eta^{\eps}  \qquad \text{when } d_{\mathcal{B}_k^c}(x) \leq 2\eta\rho^k.
    \end{equation}
Finally, by definition of $\mathcal{B}_0$ it can be seen that $d_{\mathcal{B}_0^c}(x) \geq \frac{1}{4}|x|$ whenever $|x|> 1$, so we also have
    \begin{equation}\label{major}
        \updelta_{\mathcal{B}_0^c}^{\eps}(x) \geq \frac{|x|^{\eps}}{4^\eps} \qquad \text{for }  |x|>1.
    \end{equation}
We are now ready to prove Theorem \ref{th}.
\begin{proof}
    The statement is invariant under dilation, so we assume $r_0=1$ without loss of generality. We set $\Tilde{u}:= u \chi_{B_3}$ and $L\Tilde{u}=:\Tilde{f}$ in $\Omega \cap B_2$, so that $\Tilde{u}\in L^{\infty}(\R^n)$ and $\Tilde{u}=u$ in $B_2$. By definition of $\Tilde{u}$ and Lemma \ref{wlog}, for some $C= C(n,s, \Lambda)>0$ we get the estimate
    \[ \| \Tilde{f}\|_{L^{\infty}(\Omega \cap B_2)} \leq \| f\|_{L^{\infty}(\Omega \cap B_2)} + \| L \left( u \chi_{B_3^c} \right) \|_{L^{\infty}(B_2)} \leq \| f\|_{L^{\infty}(\Omega \cap B_2)} + C\| u \|_{L^1_{2s}(\R^n)}. \]
    Hence, we may assume that 
    $\|\Tilde{f}\|_{L^{\infty}(\Omega \cap B_2)} \leq 1$ and $\|\Tilde{u}\|_{L^{\infty}(\R^n)} \leq 1$ up to dividing by a constant. Notice also that $\tilde{u} \in C(\Omega \cap B_2)$ by the interior regularity estimates of Theorem \ref{interiorestimates}. We begin by proving that for all $x_0 \in \partial \Omega \cap B_1$ there holds
    \begin{equation}\label{Holdercont}
        |\tilde{u}(x)| = |\tilde{u}(x)-\tilde{u}(x_0)| \leq |x-x_0|^{\sigma} \qquad \text{ for any } x \in B_1(x_0).
    \end{equation}
    Up to translation, we denote $x_0=0$ only for this part of the proof. To show \eqref{Holdercont} it is enough to find $0 < \eta (n, s, \lambda, \Lambda) < 1/2$, $0 < \rho (\eta) <1$ and $\sigma>0$ such that
    \begin{equation}\label{induction}
        \| \Tilde{u} \|_{L^{\infty}(\Omega \, \cap B_{\rho^k})} \leq \rho^{k \sigma} \qquad \text{for all } k \in \N.
    \end{equation}
    We prove \eqref{induction} by induction. First, when $k = 0$ the statement holds for any choice of $0 < \rho <1$ and $\sigma>0$. Then, we fix $k\in \N$ and we assume that
    \begin{equation}\label{inductionhypo}
    \| \Tilde{u} \|_{L^{\infty}(\Omega \, \cap B_{\rho^j})} \leq \rho^{j \sigma} \qquad \text{for all } j \in \N \text{ with } j \leq k.
    \end{equation}
    We show that
    \begin{equation}\label{conclusion}
    \| \Tilde{u} \|_{L^{\infty}(\Omega \, \cap B_{\rho^{k+1}})} \leq \rho^{(k+1) \sigma}. \end{equation}
    To do so, we work on the set $\Omega \, \cap B_{\rho^k}$ and we consider the associated hyperplane $H_{0, \rho^k}$ from the definition of Reifenberg flatness. We fix a coordinate system for which 
    \[ H_{0, \rho^k} = \{x \in \R^n : x_n = 0\} \]
    and we consider the barrier $v_k$ defined above in these new coordinates.
    Now, we show that the following system holds pointwise
    \begin{equation}\label{compsystem}
        \begin{cases}
            L v_k \geq 1 & \text{on } \Omega  \cap B_{\rho^k} \\ v_k \geq \tilde{u} & \text{on } \lt(\Omega  \cap B_{\rho^k} \rt)^c.
    \end{cases}
    \end{equation}
    For all $x \in \Omega  \cap B_{\rho^k}$ we have that $d_{\mathcal{B}_k^c}(x) \leq 2 \rho^k$, so $d_{\mathcal{B}_k^c}^{\eps - 2s}(x) \geq (2 \rho^k) ^{\eps -2s}$ for $\eps$ small. We impose $\sigma < \eps/2$, so from \eqref{barrierk} we obtain
    \[ L v_k(x) \geq c\, 4^{\eps}\rho^{(k-1)\sigma} d_{\mathcal{B}_k^c}^{\eps - 2s}(x) \geq c \rho^{k(\sigma+\eps-2s)} \geq (c \rho ^k)^{2\eps - 2s} \qquad \text{for all } x \in \Omega \cap B_{\rho^k}. \]
    First we choose $\eps$ such that $2\eps - 2s < 0$, then $\rho<1$ small enough such that $c \rho < 1$. In this way we get $ L v_k(x) \geq 1$ for all $x \in \Omega  \cap B_{\rho^k}$ and for all $k \in \N$, thus finding the first line of \eqref{compsystem}. Next, we check that $v_k \geq {\tilde{u}}$ over the set $\lt(\Omega \cap B_{\rho^k} \rt)^c$. We get it right away when $x \in B_1 \setminus \Omega$, because $\tilde{u}(x)=0$ and $v_k \geq 0$ on $\R^n$. If $x \in B_1^c$ then there holds $|x| > \rho^k$, so by \eqref{major} we have 
    \[ v_k(x) = 4^{\eps}\rho^{(k-1)\sigma}\updelta_{\mathcal{B}_0^c}^{\eps}\lt( \frac{x}{\rho^k} \rt) \geq  \rho^{(k-1)\sigma} \frac{|x|^{\eps}}{\rho^{k\eps}} \geq \rho^{k(\sigma - \eps)} \geq 1.\]
    Since $1 \geq \|\Tilde{u}\|_{L^{\infty}(\R^n)}$, we also get $v_k(x) \geq \tilde{u}(x)$ for $x \in B_1^c$. Finally, if $x \in \Omega \cap (B_1 \setminus B_{\rho^k})$ then
    \[ |x| \in [ \rho^j, \rho^{j-1}) \qquad \text{for some } j\in \N \text{ with } 1 \leq j \leq k. \]
    By induction hypothesis \eqref{inductionhypo}, there holds $ \tilde{u}(x) \leq \rho^{(j-1)\sigma}$. Since $|x|> \rho^k$, again by \eqref{major} we obtain 
    \[ v_k(x) \geq \rho^{(k-1)\sigma} \frac{|x|^{\eps}}{\rho^{k\eps}} \geq \rho^{(k-1)\sigma + (j-k)\eps}.\] 
    Altogether, we have
    \[ \frac{v_k(x)}{\tilde{u}(x)} \geq \rho^{(j-k)(\eps-\sigma)} \geq 1 \qquad \text{since} \quad j-k \leq 0,\quad \eps-\sigma >0\, \quad\text{and}\quad 0<\rho <1.\]
    Thus we deduce $v_k \geq \tilde{u}$ on $\lt(\Omega \cap B_{\rho^k} \rt)^c$ as well as the validity of the second line of \eqref{compsystem}. Now, we aim at applying the weak comparison principle to get $\tilde{u} \leq v_k$ also on $\Omega \cap B_{\rho^k}$: in order to do so we must show that $L v_k \geq L\tilde{u}$ weakly on $\Omega \cap B_{\rho^k}$. We already know that $1 \geq L \tilde{u}$ in the weak sense and that $ L v_k(x) \geq 1$ pointwise over $\Omega \cap B_{\rho^k}$, so it is only left to prove that 
    \[ \langle v_k, v_k \rangle_{K, \Omega \cap B_{\rho^k}} < + \infty.\]
    However, we cannot argue directly because we have no information on the behavior of $\partial \Omega$ inside the stripe $\lt \{ x \in B_{\rho^k}: d(x, H_{0,\rho^k}) \leq \eta \rho^k \rt \}$. More precisely, we do not know whether $d(\mathcal{B}_k, \Omega \cap B_{\rho^k})$ is strictly positive. To solve this, it is enough to consider the set
    \[ \mathcal{B}_{0,t} := \overline{B \lt(- \lt(\frac{1}{8} +\eta \rt) e_n, \frac{1}{8}-t \rt)} \qquad \text{for small } t>0.\]
    Setting $\mathcal{B}_{k,t} := \rho^k \mathcal{B}_{0,t}$, we apply Lemma \ref{finalbarrier} as before and we define the barriers $v_{k,t}$ for all $k \in \N$ accordingly. Since $d(\mathcal{B}_{k,t}, \Omega \cap B_{\rho^k})>0$, now we have $v_{k,t}\in C^{\infty}(\mathcal{B}_{k,t}^c) \cap \text{Lip}(\mathcal{B}_{k,t/2}^c)$ for $k\in\N$. Hence, denoting $\Omega_k := \Omega \cap B_{\rho_k}$ for short:
    \[ \frac{1}{\Lambda}\langle v_{k,t}, v_{k,t} \rangle_{K, \Omega_k} \leq \int_{\Omega_k \times \Omega_k} \frac{(v_{k,t}(x)-v_{k,t}(y))^2}{|x-y|^{n+2s}}dxdy + 2\int_{\Omega_k \times \Omega_k^c} \frac{(v_{k,t}(x)-v_{k,t}(y))^2}{|x-y|^{n+2s}}dxdy = I_1 + 2 I_2. \]
    The function $ v_{k,t}$ is Lipschitz over $\Omega_k$, so we get
    \begin{equation}\label{energy}
        I_1 \leq C \int_{\Omega_k \times \Omega_k} \frac{1}{|x-y|^{n+2s-2}}\,dxdy < + \infty.
    \end{equation}
    Next, we introduce the set $A_k := B_{2 \rho_k} \setminus \{ x \in \R^n : x_n \leq -(\eta + \frac{t}{2})\rho^k\}$ and we write
    \[ I_2 = \int_{\Omega_k \times (\Omega_k^c \cap A_k)} \frac{(v_{k,t}(x)-v_{k,t}(y))^2}{|x-y|^{n+2s}}dxdy + \int_{\Omega_k \times (\Omega_k^c \cap A_k^c)} \frac{(v_{k,t}(x)-v_{k,t}(y))^2}{|x-y|^{n+2s}}dxdy = I_3 + I_4. \]
    The function $v_{k,t}$ is Lipschitz also over $A_k$, so we argue as in \eqref{energy} and obtain $I_3 < + \infty$. Now, we majorize $I_4$:
    \[ I_4 \leq 2 \int_{\Omega_k \times A_k^c} \frac{v_{k,t}(x)^2}{|x-y|^{n+2s}}dxdy + 2 \int_{\Omega_k \times A_k^c} \frac{v_{k,t}(y)^2}{|x-y|^{n+2s}}dxdy = 2I_5 + 2I_6. \]
    For $x \in \Omega_k$ we have $v_{k,t}(x) \leq \| v_{k,t}\|_{L^{\infty}(\Omega_k)} < + \infty$, thus we find $I_5 < +\infty$. Finally, by construction it is possible to see that for $C$ depending on $\sigma$, $\eps$, $\rho$ and $k$ (which here are all fixed) there holds
    \begin{equation}\label{lastestimate}
        v_{k,t}(y) \leq C d^{\eps}(y, \mathcal{B}_{k,t}) \leq C|x-y|^{\eps} \qquad \text{for all } x \in \Omega_k \text{ and } |y| \geq 3\rho^k.
    \end{equation}
    Thus,
    \[ I_6 = \int_{\Omega_k \times (A_k^c \cap B_{3\rho^k})} \frac{v_{k,t}(y)^2}{|x-y|^{n+2s}}dxdy + \int_{\Omega_k \times B^c_{3\rho^k}} \frac{v_{k,t}(y)^2}{|x-y|^{n+2s}}dxdy = I_7 + I_8. \]
    Since $v_{k,t}(y)$ is bounded over $B_{3\rho^k}$ we have that $I_7 < + \infty$. On the other hand, by \eqref{lastestimate}:
    \[ I_8 \leq C \int_{\Omega_k \times B_{3 \rho^k}^c}  \frac{1}{|x-y|^{n+2s-2\eps}}\,dxdy < + \infty. \]
    Therefore, there holds $\langle v_{k,t}, v_{k,t} \rangle_{K, \Omega \cap B_{\rho^k}} < + \infty$ for all $t>0$ small enough. Noticing that \eqref{major} is true for $v_{0,k}$ as well, we can follow the previous computations and see that the system
    \begin{equation}\label{alsodistr} \begin{cases}
            L v_{k,t} \geq L\tilde{u} & \text{on } \Omega  \cap B_{\rho^k} \\ v_{k,t} \geq \tilde{u} & \text{on } \lt(\Omega  \cap B_{\rho^k} \rt)^c
    \end{cases}\end{equation}
    holds in the weak sense. Therefore, applying Theorem \ref{maxprinciple} we deduce that $\tilde{u} \leq v_{k,t}$ on $\Omega \cap B_{\rho^k}$. Letting $t \rightarrow 0$ we conclude $\tilde{u} \leq v_k$ on $\Omega \cap B_{\rho^k}$, namely
    \[ \tilde{u}(x) \leq 4^{\eps} \rho^{(k-1)\sigma} \updelta_{\mathcal{B}_0^c}^{\eps}\lt( \frac{x}{\rho^k} \rt)\qquad \text{for all} \,\, x \in \Omega \cap B_{\rho^k}. \]
    For all $x \in \Omega \cap B_{\eta \rho^k}$ there holds $d_{\mathcal{B}_k^c}(x) \leq 2\eta\rho^k$, so we have by \eqref{holderbarrier}
    \[ \tilde{u}(x) \leq 4^{\eps} C_H \rho^{(k-1)\sigma} \eta^{\eps} \qquad \text{for all} \,\, x \in \Omega \cap B_{\eta \rho^k}. \]
    Setting $\eta=\rho$, we find
    \[ \tilde{u}(x) \leq 4^{\eps} C_H \rho^{(k-1)\sigma} \rho^{\eps} =\rho^{(k+1)\sigma} \lt( 4^{\eps} C_H \rho^{\eps - 2 \sigma} \rt) \qquad \text{for all} \,\, x \in \Omega \cap B_{\rho^{k+1}}. \]
    Recalling that $\eps - 2 \sigma > 0$, we can choose $\rho < 1/2$ small enough such that $4^{\eps} C_H \rho^{\eps - 2 \sigma} \leq 1$ and we obtain
    \[ \tilde{u}(x) \leq  \rho^{(k+1)\sigma} \qquad \text{for all} \,\, x \in \Omega \cap B_{\rho^{k+1}}. \]
    To complete the proof of \eqref{conclusion} we only need to show the lower bound $\tilde{u}(x) \geq -\rho^{(k+1)\sigma}$ for all $x \in \Omega \cap B_{\rho^{k+1}}$. The argument is the same as before, except that in this case we use the function $-v_k$ as a barrier thus finding $\tilde{u} \geq - v_k$. Hence we get \eqref{conclusion} and we conclude \eqref{induction} by induction, obtaining relation \eqref{Holdercont} as well. In particular, for all $x \in B_{3/4}$ we choose $x_0\in \partial \Omega$ such that $|x-x_0| = d_{\Omega}(x)$ and we infer
    \begin{equation}\label{pointholder}
        |u(x)|=|\tilde{u}(x)| \leq d_{\Omega}^{\sigma}(x) \qquad \text{for all} \,\, x \in B_{3/4}.
    \end{equation}
    Now, it remains to combine \eqref{pointholder} with the interior estimates of Theorem \ref{interiorestimates} to obtain the $\sigma$-H\"older regularity result. Let $x, y \in \Omega \cap B_{1/2}$ and assume that $r := d_{\Omega}(x) \leq d_{\Omega}(y)$. Setting $\zeta := |x - y|$, if $\zeta \geq \frac{r}{2}$ then
    \[ |u(x)- u(y)| \leq |u(x)|+|u(y)| \leq r^{\sigma} + (\zeta + r)^{\sigma} \leq C \, \zeta^{\sigma} = C|x - y|^{\sigma}\]
    and we immediately find $\sigma$-H\"older regularity. If instead $\zeta < \frac{r}{2}$ then we have $B_r(x) \subset \Omega \cap B_1$ and $y \in B_{r/2}(x)$. Consider the rescaled function $\tilde{u}_r(z) := \tilde{u}(x + rz)$ and the operator $L_r \in \mathcal{L}_s^n(\lambda, \Lambda)$ with kernel $K_r(y):= r^{n+2s}K(ry)$, by relation \eqref{rescaledop} we obtain
    \[ (L_r \tilde{u}_r)(z)= r^{2s}(L\tilde{u})(x+rz)= r^{2s}\tilde{f}(x+rz)=:g_r(z) \qquad \text{for } z \in B_1.\]
    Since $\| \tilde{f}\|_{L^{\infty}(\Omega \cap B_2)} \leq 1$, notice that $\|g_r\|_{L^{\infty}(B_1)}<+\infty$. Let now $\Bar{x}\in \partial \Omega$ such that $r = |x - \Bar{x}|$, thanks to estimate \eqref{pointholder} we get for any $z \in B_{1/2r}$
    \[ |\tilde{u}_r(z)| \leq d_{\Omega}^{\sigma}(x + rz) \leq (|x - \overline{x}| + r|z|)^{\sigma} \leq C (|x - \overline{x}|^{\sigma} + r^{\sigma}|z|^{\sigma}) = C r^{\sigma}(1 + |z|^{\sigma}). \]
    Given that $\| \tilde{u}_r \|_{L^{\infty}(\R^n)} = \| \tilde{u} \|_{L^{\infty}(\R^n)} \leq 1$, we deduce $\tilde{u}_r \in L^1_{2s}(\R^n)$ with $\| \tilde{u}_r \|_{L^1_{2s}(\R^n)} \leq C r^{\sigma}$.
    Altogether we showed that $\tilde{u}_r \in L^1_{2s}(\R^n)$ is a distributional solution to the equation $L_r \tilde{u}_r = g_r$ in $B_1$, so by Theorem \ref{interiorestimates} with $C^{\sigma}$ norm on the left-hand side we have the estimate
    \[ \|\tilde{u}_r \|_{C^{\sigma}(B_{1/2})} \leq C \lt( \| \tilde{u}_r \|_{L^1_{2s}(\R^n)} + \,\| g_r \|_{L^{\infty}(B_1)} \rt) \leq C r^{\sigma}. \]
    Since $\|\tilde{u}_r \|_{C^{\sigma}(B_{1/2})} = r^{\sigma} \|\tilde{u} \|_{C^{\sigma}(B_{r/2}(x))}$, we obtain $\|\tilde{u} \|_{C^{\sigma}(B_{r/2}(x))} \leq C$ and $|\tilde{u}(x) - \tilde{u}(y)| \leq C r^{\sigma}$. Recalling that $u = \tilde{u}$ in $B_{1/2}$, for any $x, y \in B_{1/2}$ there holds
    \[ |u(x)- u(y)| \leq C|x - y|^{\sigma}, \]
    which concludes the proof.
\end{proof}
\begin{remark}\label{rem}
    The previous theorem holds for continuous distributional solutions as well. Indeed, once \eqref{alsodistr} is true in the weak sense then it is enough to apply the maximum principle for distributional solutions \cite[Lemma 2.3.5]{fernandez2024integro} and to proceed in the same way. More in general, we proved that weak solutions to nonlocal equations over Reifenberg flat sets are continuous up to the boundary, hence from now on we may work with distributional solutions which are continuous up to the boundary without loss of generality.
\end{remark}

We now show $C^{s-\eps}$ regularity in the same framework of Theorem \ref{th}. The proof goes by compactness and it is an immediate adaptation to the Reifenberg flat case of the argument of \cite[Theorem 6.8]{ros2024optimal}, where the same conclusion is drawn for $\delta$-Lipschitz domain with sufficiently small $\delta$. We highlight that Theorem \ref{compactness} is sharp in view of the optimality of \cite[Theorem 6.8]{ros2024optimal}.

\begin{theorem}\label{compactness}
    Let $ L \in \mathcal{L}_s^n(\lambda, \Lambda)$ and $\gamma > 0$. Let $\Omega \subset \R^n$ be an open and connected $(\eta, r_0)$-Reifenberg flat set for some $r_0>0$, $f \in L^{\infty}(\Omega \cap B_1)$ and $u \in C(B_1) \cap L^1_{2s}(\R^n)$ be a distributional solution to
    \[ \begin{cases}
        Lu = f & \text{on } \Omega \cap B_1 \\ 
        u=0 & \text{on } B_1 \setminus \Omega.
    \end{cases} \]
    Then, there exist $\eta_0 = \eta_0(n,s,\lambda, \Lambda, \gamma)>0$ and  $C= C(n,s, \lambda, \Lambda, \gamma)>0$ such that for all $\eta \leq \eta_0$ we have $u \in C^{s-\gamma}_{loc}(B_1)$ with
   \begin{equation}\label{localweakbound}
   \|u\|_{C^{s - \gamma}(B_{1/2})} \leq C \lt( \| u \|_{L^1_{2s}(\R^n)} + \,\| f\|_{L^{\infty}(\Omega \cap B_1)} \rt). \end{equation}
\end{theorem}

\begin{proof}
    We set $r_0=1$ without loss of generality. Applying the same renormalization argument at the beginning of the proof of Theorem \ref{th}, we assume $u \in L^{\infty}(\R^n)$ with $\|u\|_{L^{\infty}(\R^n)} \leq 1$ and $\|f\|_{L^{\infty}(\Omega \cap B_1)} \leq 1$. Our goal is to show that there exists $C = C(n,s, \gamma, \lambda, \Lambda)>0$ such that for any $x_0 \in \partial \Omega \cap B_{1/2}$ we have
    \begin{equation}\label{contradiction}
        |u(x)| \leq C |x-x_0|^{s - \gamma} \qquad \text{for all} \,\, x \in B_{1/2} \cap \Omega.
    \end{equation}
    We assume by contradiction that \eqref{contradiction} does not hold. Therefore, there exists sequences $\{\eta_k\}_{k \in \N}$ with $\eta_k \rightarrow 0$ as $k \rightarrow + \infty$, $\{\Omega^k\}_{k \in \N}$ $(\eta_k,1)$-Reifenberg flat sets with $0 \in \partial \Omega^k$, $\{L_k\}_{k \in \N} \subset \mathcal{L}_s^n(\lambda, \Lambda)$ with kernels $\{K_k\}_{k \in \N}$, $\{f_k\}_{k \in \N} \subset L^{\infty}(\Omega^k \cap B_1)$ and  $\{u_k\}_{k \in \N} \subset C(B_1) \cap L^{\infty}(\R^n)$ with
    \begin{equation}\label{normalization}
        \|u_k\|_{L^{\infty}(\R^n)} + \|f_k\|_{L^{\infty}(\Omega \cap B_1)} \leq 1
    \end{equation}
    and such that $u_k$ is a distributional solution to
    \[ \begin{cases}
        L_k u_k = f_k & \text{on } \Omega^k \cap B_1 \\ 
        u_k =0 & \text{on } B_1 \setminus \Omega^k,
    \end{cases}\]
    but we also have
    \[ \sup_{k \in \N} \sup_{r>0} \frac{\| u_k \|_{L^{\infty}(B_{r/2})}}{r^{s - \gamma}} = + \infty. \]
    By \cite[Lemma 4.4.11]{fernandez2024integro}, up to extraction of a subsequence there exists $\{\eta_k\}_{k \in \N}$ with $\eta_k \rightarrow 0$ as $k \rightarrow + \infty$ and such that
    \begin{equation}\label{absurdest}
        \| u_k \|_{L^{\infty}(B_{r_k})} \geq C r^{s-\gamma}
    \end{equation}
    for a constant $C>0$ independent of $k$. Considering the rescaled functions
    \[ v_k(x):= \frac{u_k(r_k x)}{\| u_k \|_{L^{\infty}(B_{r_k})}} \qquad \text{and} \qquad \tilde{f}_k(x) := \frac{r^{2s}_k f_k(r_k x)}{\| u_k \|_{L^{\infty}(B_{r_k})}},  \]
    we notice that
    \begin{equation}\label{L^1_2sconv}
        \| v_k \|_{L^{\infty}(B_1)}=1 \qquad \text{and} \qquad |v_k(x)| \leq 2(1 + |x|^{s-\gamma}) \quad \text{for all} \,\, x \in \R^n.
    \end{equation} 
    Let $\{\tilde{L}_k\}_{k \in \N} \subset \mathcal{L}_s^n(\lambda, \Lambda)$ be operators with kernels $\tilde{K}_k(y) := r_k^{n+2s} K_k(r_k y)$. Then, by the scaling properties of the class $\mathcal{L}_s^n(\lambda, \Lambda)$, $v_k$ is a distributional solution to
    \[ \begin{cases}
        \tilde{L}_k v_k = \tilde{f}_k & \text{on } 
        \Omega^k_{1/r^k} \cap B_{1/r^k} \\ 
        v_k =0 & \text{on } B_{1/r^k} \setminus \Omega^k_{1/r^k}.
    \end{cases}\]
    Furthermore, by \eqref{normalization} and \eqref{absurdest} we have
    \begin{equation}\label{vanish}
        \| \tilde{f}_k \|_{L^{\infty} \lt( \Omega^k_{1/r^k} \cap B_{1/r^k} \rt)} \leq \| f_k \|_{L^{\infty} \lt( \Omega^k \cap B_1 \rt)} C^{-1} r_k^{s + \gamma} \longrightarrow 0 \qquad \text{as} \qquad k \rightarrow + \infty.
    \end{equation}
    Now, Theorem \ref{th} together with \eqref{normalization} and \eqref{vanish} implies that that the sequence $\{v_k\}_{k \in \N}$ is uniformly bounded in $C^{\sigma}( \Omega^k_{1/r^k} \cap B_{1/2r^k})$. By Ascoli-Arzelà theorem, it converges locally uniformly to some $v \in C(\R^n) \cap L^1_{2s}(\R^n)$ satisfying 
    \begin{equation}\label{absurd} \| v\|_{L^{\infty}(B_1)}=1 \qquad \text{and} \qquad |v(x)| \leq 2 (1 + |x|^{s- \gamma}) \quad \text{for all}\,\, x \in \R^n.\end{equation}
    Since $\eta_k, r_k \rightarrow 0$, up to another extraction there exists $e \in \S^{n-1}$ such that $\Omega^k_{1/r^k} \cap B_{1/r^k} \rightarrow \{ x \cdot e >0\}$ as $k \rightarrow + \infty$. In addition, respectively from \eqref{L^1_2sconv} and \eqref{vanish} we get $v_k \longrightarrow v$ and $\tilde{f}_k \longrightarrow 0$ in $L^1_{2s}(\R^n)$ as $k \rightarrow + \infty$. Thanks to the stability of distributional solutions \cite[Proposition 2.2.36]{fernandez2024integro}, there exists $L_{\infty} \in \mathcal{L}_s^n(\lambda, \Lambda)$ such that the following system holds in distributional sense
    \[ \begin{cases}
        L_{\infty} v = 0 & \text{on } 
        \{ x \cdot e >0\} \\ 
        v =0 & \text{on } \{ x \cdot e \leq 0\}.
    \end{cases}\]
    By Liouville theorem in the half-space \cite[Theorem 6.2]{ros2024optimal} and the second property in \eqref{absurd} we have that $v=0$, but this is a contradiction with $\|v\|_{L^{\infty}(B_1)}=1$ so \eqref{contradiction} must be true. We complete the proof combining \eqref{contradiction} with the interior regularity estimates, arguing as in Theorem \ref{th}.
\end{proof}
In the end, the proof of Theorem \ref{goal} is a direct consequence of the previous result.
\begin{proof}[Proof of Theorem \ref{goal}]
    Since $\Omega$ is bounded, we cover it with finitely many balls of radius $1/2$ and we apply Theorem \ref{compactness} to each of them, thus obtaining $u\in C^{s-\eps}(\overline{\Omega})$. Estimate \eqref{weakbound} comes directly from its local counterpart \eqref{localweakbound} and the $L^{\infty}$ bound for weak solutions of \cite[Lemma 2.3.9]{fernandez2024integro}.
\end{proof}
\begin{remark}
    In view of Remark \ref{rem}, Theorem \ref{goal} is true also for distributional solutions which are continuous up to the boundary.
\end{remark}
\end{section}

{\bf Acknowledgements.} We kindly thank Michael Goldman for recommending this problem and for the insightful discussions leading to the result. Part of this work was carried out during a visit at University of Pisa which was supported by a public grant from the Fondation Math\'ematique Jacques Hadamard and by the project G24-202 ``Variational methods for geometric and optimal matching problems'' funded by Università Italo Francese.

\bibliographystyle{acm}
\bibliography{prep}

\end{document}